\documentclass[12pt]{article}
\usepackage{myart}
\usepackage{amsfonts}

\normalbaselines
\input{tcilatex}

\begin{document}

\title{Deformation principle and problem of parallelism in geometry and physics.}
\author{Yuri A.Rylov}
\date{Institute for Problems in Mechanics, Russian Academy of Sciences \\
101-1 ,Vernadskii Ave., Moscow, 119526, Russia \\
email: rylov@ipmnet.ru\\
Web site: {$http://rsfq1.physics.sunysb.edu/\symbol{126}rylov/yrylov.htm$}}
\maketitle

\begin{abstract}
The deformation principle admits one to obtain a very broad class of
nonuniform geometries as a result of deformation of the proper Euclidean
geometry. The Riemannian geometry is also obtained by means of a deformation
of the Euclidean geometry. Application of the deformation principle appears
to be not consecutive, and the Riemannian geometry appears to be not
completely consistent. Two different definitions of two vectors parallelism
are investigated and compared. The first definitions is based on the
deformation principle. The second definition is the conventional definition
of parallelism, which is used in the Riemannian geometry. It is shown, that
the second definition is inconsistent. It leads to absence of absolute
parallelism in Riemannian geometry and to discrimination of outcome outside
the framework of the Riemannian geometry at description of the space-time
geometry.
\end{abstract}


PACS2003: 02.40Ky


\newpage

\section{Introduction}

There are two different approaches to geometry: mathematical and physical
ones. In the mathematical approach a geometry is a construction founded on a
system of axioms about points and straight lines. Practically any system of
axioms, containing concepts of a point and a straight line, may be called a
geometry. Well known mathematician Felix Klein \cite{K37} supposed that only
such a construction on a point set is a geometry, where all points of the
set have the same properties (uniform geometry). For instance, Felix Klein
insisted that Euclidean geometry and Lobachevsky geometry are geometries,
because they are uniform, whereas the Riemannian geometries are not
geometries at all. As a rule the Riemannian geometries are not uniform, and
their points have different properties. According to the Felix Klein
viewpoint, they should be called as ''Riemannian topographies'' or as
''Riemannian geographies''. Thus, at the mathematical approach to geometry
the main feature of geometry is existence of some axiomatics. One may say
that the mathematical geometry (mathematical approach to geometry) is a
system of axioms. Practically one can construct axiomatics only for uniform
geometries, and any mathematical geometry is a uniform geometry.

Riemannian geometries are not uniform geometries, in general. Practically
one cannot construct axiomatics for each of the continuous set of Riemannian
geometries, and any Riemannian geometry is obtained as a result of some
deformation of the proper Euclidean geometry, when the infinitesimal
Euclidean interval $ds_{{\rm E}}^{2}$ is replaced by the infinitesimal
Riemannian interval $ds^{2}=g_{ik}dx^{i}dx^{k}$. Such a replacement is a
deformation of the Euclidean space.

Such an approach to geometry, when a geometry is obtained as a result of
deformation of the proper Euclidean geometry will be referred to as the
physical approach to geometry. The obtained geometry will be referred to as
physical geometry. The physical geometry has not its own axiomatics. It uses
''deformed '' Euclidean axiomatics. The term ''physical geometry'' is used,
because it is very convenient for application to physics and can be used as
a space-time geometry. Felix Klein referred to a physical geometry as a
topography, but we think that another name is important not in itself, but
only because it describes another method of the geometry construction.

Physical geometry describes mutual disposition of geometric objects in the
space, or mutual dispositions of events in the event space (space-time). The
mutual dispositions is described by the distance between any two points. It
is of no importance, whether the geometry has any axiomatics or not. One may
say that the physical geometry (physical approach to geometry) is a
conception, describing mutual dispositions of geometric objects and points.
Physical geometry may be not uniform, and it is not uniform in many cases.
Metric $\rho $ (distance between two points) is a unique characteristic of a
physical geometry. World function $\sigma =\frac{1}{2}\rho ^{2}$ \cite{S60}
is more convenient for description of a physical geometry, because it is
real even for the space-time, where $\rho =\sqrt{2\sigma }$ may be
imaginary. Besides, usually the term metric is associated with some
constraints on metric (triangle axiom, positivity of $\rho $). The term
world function does not associate with these constraints directly.

Attempts of construction of a metric geometry without the constraint, imosed
by the triangle axiom, were made earlier \cite{M28}. It is so called
distance geometry \cite{B53}. Unfortunately, these attempts did not lead to
a construction of a pithy geometry in terms of only metric.

Construction of any physical geometry is determined by {\it the deformation
principle}. It works as follows. The proper Euclidean geometry ${\cal G}_{%
{\rm E}}$ can be described in terms and only in terms of the world function $%
\sigma _{{\rm E}}$, provided $\sigma _{{\rm E}}$ satisfies some constraints
formulated in terms of $\sigma _{{\rm E}}$ \cite{R02}. It means that all
geometric objects ${\cal O}_{{\rm E}}$ can be described $\sigma $-immanently
(i.e. in terms of $\sigma _{{\rm E}}$ and only of $\sigma _{{\rm E}}$) $%
{\cal O}_{{\rm E}}={\cal O}_{{\rm E}}\left( \sigma _{{\rm E}}\right) $.
Relations between geometric objects are described by some expressions ${\cal %
R}_{{\rm E}}={\cal R}_{{\rm E}}\left( \sigma _{{\rm E}}\right) $. Any
physical geometry ${\cal G}_{{\rm A}}$ can be obtained from the proper
Euclidean geometry by means of deformation, when the Euclidean world
function $\sigma _{{\rm E}}$ is replaced by some other world function $%
\sigma _{{\rm A}}$ in all definitions of Euclidean geometric objects ${\cal O%
}_{{\rm E}}={\cal O}_{{\rm E}}\left( \sigma _{{\rm E}}\right) $ and in all
Euclidean relations ${\cal R}_{{\rm E}}={\cal R}_{{\rm E}}\left( \sigma _{%
{\rm E}}\right) $ between them. As a result we have the following change 
\[
{\cal O}_{{\rm E}}={\cal O}_{{\rm E}}\left( \sigma _{{\rm E}}\right)
\rightarrow {\cal O}_{{\rm A}}={\cal O}_{{\rm E}}\left( \sigma _{{\rm A}%
}\right) ,\qquad {\cal R}_{{\rm E}}={\cal R}_{{\rm E}}\left( \sigma _{{\rm E}%
}\right) \rightarrow {\cal R}_{{\rm A}}={\cal R}_{{\rm E}}\left( \sigma _{%
{\rm A}}\right) 
\]
The set of all geometric objects ${\cal O}_{{\rm A}}$ and all relations $%
{\cal R}_{{\rm A}}$ between them forms a physical geometry, described by the
world function $\sigma _{{\rm A}}$. Index ${\rm E}$ in the relations of
physical geometry ${\cal G}_{{\rm A}}$ means that axiomatics of the proper
Euclidean geometry was used for construction of geometric objects ${\cal O}_{%
{\rm E}}={\cal O}_{{\rm E}}\left( \sigma _{{\rm E}}\right) $ and of
relations between them ${\cal R}_{{\rm E}}={\cal R}_{{\rm E}}\left( \sigma _{%
{\rm E}}\right) $. The same axiomatics is used for all geometric objects $%
{\cal O}_{{\rm A}}={\cal O}_{{\rm E}}\left( \sigma _{{\rm A}}\right) $ and
relations between them ${\cal R}_{{\rm A}}={\cal R}_{{\rm E}}\left( \sigma _{%
{\rm A}}\right) $ in the geometry ${\cal G}_{{\rm A}}$. But now this
axiomatics has another form, because of deformation $\sigma _{{\rm E}%
}\rightarrow \sigma _{{\rm A}}$. It means that the proper Euclidean geometry 
${\cal G}_{{\rm E}}$ is the basic geometry for all physical geometries $%
{\cal G}$ obtained by means of a deformation of the proper Euclidean
geometry. If basic geometry is fixed (it is this case that will be
considered further), the geometry on the arbitrary set $\Omega $ of points
is called T-geometry (tubular geometry). The T-geometry is determined \cite
{R02} by setting the world function $\sigma $: 
\begin{equation}
\sigma :\;\;\;\Omega \times \Omega \rightarrow {\Bbb R},\qquad \sigma \left(
P,P\right) =0,\qquad \forall P\in \Omega  \label{a1}
\end{equation}
In general, no other constraints are imposed, although one can impose any
additional constraints to obtain a special class of T-geometries. T-geometry
is symmetric, if in addition 
\begin{equation}
\sigma \left( P,Q\right) =\sigma \left( Q,P\right) ,\qquad \forall P,Q\in
\Omega  \label{a1.1}
\end{equation}

Deformation ${\Bbb R}^{n}\rightarrow \Omega $ of the $n$-dimensional proper
Euclidean space to an arbitrary set $\Omega $ of points is a deformation in
the broad sense. This deformation can change the dimension of a geometric
object and the dimension of the whole space. For instance, the resulting
T-geometry does not depend on the dimension $n$ of deformed proper Euclidean
space. Only final world function $\sigma $ is important for the T-geometry
properties. This admits one to consider T-geometry as something
self-sufficient and to ignore the deformation which produces T-geometry from
the Euclidean geometry.

The Riemannian geometry is a physical geometry. It is constructed on the
basis of the {\it deformation principle}, i.e. in the same way as
T-geometry. But class of possible Riemannian deformations is not so general
as the class of all possible deformations. It is restricted by the
constraint 
\begin{equation}
\sigma _{{\rm R}}\left( x,x^{\prime }\right) =\frac{1}{2}\left( \int\limits_{%
{\cal L}_{\left[ xx^{\prime }\right] }}\sqrt{g_{ik}dx^{i}dx^{k}}\right) ^{2}
\label{c1}
\end{equation}
where $\sigma _{{\rm R}}$ is the world function of Riemannian geometry, and $%
{\cal L}_{\left[ xx^{\prime }\right] }$ denotes segment of geodesic between
the points $x$ and $x^{\prime }$. The Riemannian geometry is determined by
the dimension $n$ and $n\left( n+1\right) /2$ functions $g_{ik}$ of one
point $x$, whereas the class of possible T-geometries is determined by one
function $\sigma $ of two points $x$ and $x^{\prime }$.

A use of the deformation principle is sufficient for a construction of any
physical geometry. All relations between geometric objects appear to be as
consistent as they are consistent in the proper Euclidean geometry. The
deformation principle does not use any logical conclusions and leads to a
construction of a consistent physical geometry. Moreover, a use of
additional means of the geometry construction is undesirable, because these
means may disagree with the deformation principle. In the case of such a
disagreement the obtained geometry appears to be inconsistent.

Although the Riemannian geometry is a kind of physical geometry, at its
construction one uses additional means of description (dimension, concept of
a curve, coordinate system, continuous manifold). Some of them appear to be
incompatible with the principle of the geometry deformation, and as a result
the Riemannian geometry appears to be inconsistent. Constraint (\ref{c1}),
imposed on the world function of Riemannian geometry, restricts the class of
possible physical geometries and reduces this inconsistency, but it fails to
eliminate inconsistency completely. The $\sigma $-Riemannian geometry, i.e.
the physical geometry, constructed by means of only the deformation
principle on $n$-dimensional manifold and restricted by the constraint \cite
{S60} 
\begin{equation}
\sigma _{i}\left( x,x^{\prime }\right) g^{ik}\left( x\right) \sigma
_{k}\left( x,x^{\prime }\right) =2\sigma \left( x,x^{\prime }\right) ,\qquad
\sigma _{i}\left( x,x^{\prime }\right) \equiv \frac{\partial \sigma \left(
x,x^{\prime }\right) }{\partial x^{i}}  \label{c2}
\end{equation}
which is equivalent to (\ref{c1}), is rather close to the Riemannian
geometry. Nevertheless, the absolute parallelism is absent in the Riemannian
geometry, but it takes place in the $\sigma $-Riemannian geometry. This
difference means that the Riemannian geometry is inconsistent, because the $%
\sigma $-Riemannian geometry cannot be inconsistent.

From viewpoint of the deformation principle this inconsistency is
conditioned by a use of special properties of the world function $\sigma _{%
{\rm E}}$ of $n$-dimensional proper Euclidean space. It means as follows.
Before deformation the geometric objects ${\cal O}_{{\rm E}}$ and the
relations ${\cal R}_{{\rm E}}$ of the proper Euclidean geometry are to be
represented in the $\sigma $-immanent form. Representing ${\cal O}_{{\rm E}}$
and ${\cal R}_{{\rm E}}$ in terms of $\sigma _{{\rm E}}$, we must not use
special properties of Euclidean world function $\sigma _{{\rm E}}$. These
special properties of $\sigma _{{\rm E}}$ are formulated for $n$-dimensional
Euclidean space and contain a reference to the space dimension $n$. If these
properties are used at the description of ${\cal O}_{{\rm E}}$, or ${\cal R}%
_{{\rm E}}$, the description contains a reference to the dimension $n$ of
the space. In this case after deformation we attribute some properties of $n$%
{\it -dimensional} proper Euclidean geometry to the constructed physical
geometry. Formal criterion of application of special properties of $\sigma _{%
{\rm E}}$ is a reference to the dimension $n$. Being transformed to $\sigma $%
-immanent form, such a description of ${\cal O}_{{\rm E}}$, or ${\cal R}_{%
{\rm E}}$ contains additional points which are not characteristic for ${\cal %
O}_{{\rm E}}$, or ${\cal R}_{{\rm E}}$. Practically, these additional points
describe the coordinate system, and number of these points depends on the
space dimension $n$.

Inconsistency of the Riemannian geometry manifests itself in the parallelism
problem. The definition of $\ $two vectors parallelism in Riemannian
geometry has two defects:

\begin{enumerate}
\item  Definition of parallelism in Riemannian geometry is coordinate
dependent, because it contains a reference to the number of coordinates
(space dimension).

\item  Parallelism is defined only for two infinitesimally close vectors.
Parallelism of two remote vectors at points $P_{1}$ and $P_{2}$ is defined
by means of a parallel transport along some curve connecting points $P_{1}$
and $P_{2}$. In curved space the result of parallel transport depends on the
path of transport, and the absolute parallelism is absent, in general.
\end{enumerate}

The problem of definition of two vectors parallelism is very important,
because parallelism lies in the foundation of the particle dynamics. For
instance, in the curved space-time the free particle motion is described by
the geodesic equation 
\begin{equation}
d\dot{x}^{i}=-\Gamma _{kl}^{i}\dot{x}^{k}dx^{l},\qquad dx^{l}=\dot{x}%
^{i}d\tau  \label{a0}
\end{equation}
where $\Gamma _{kl}^{i}$ is the Christoffel symbol. Equations (\ref{a0})
describe parallel transport of the particle velocity vector $\dot{x}^{i}$
along the direction $dx^{i}=\dot{x}^{i}d\tau $ determined by the velocity
vector $\dot{x}^{i}$. If the parallel transport (\ref{a0}) appears to be
incorrect and needs a modification, the equation of motion of a free
particle needs a modification also. For instance, if at the point $x$ a set
of vectors $u^{i}$, which are parallel to the velocity vector $\dot{x}^{i}$,
appears to be consisting of many mutually noncollinear vectors, the parallel
transport of the velocity vector $\dot{x}^{i}$ stops to be single-valued,
and the world line of a free particle becomes to be random.

Definition of the scalar product of two vectors in Riemannian geometry
contains special properties of Euclidean world function and attributes to
Riemannian geometry some properties of the proper Euclidean geometry, mainly
one-dimensionality of straight lines (geodesics). In general, this
definition of scalar product is incompatible with the deformation principle.
Restriction\ (\ref{c1}), imposed on the world function $\sigma _{{\rm R}}$,
eliminates geometries admitting non-one-dimensional ''straight lines'' and
eliminates some corollaries of this incompatibility, but not all. Creators
of the Riemannian geometry tried to conserve one-dimensional straight lines
(geodesics) in the Riemannian geometry. They had achieved this goal, but not
completely, because only straight lines (geodesic) ${\cal L}\left( P_{0},%
{\bf P}_{0}{\bf P}_{1}\right) $, drawn through the point $P_{0}$ parallel to
the vector ${\bf P}_{0}{\bf P}_{1}$, is one-dimensional, whereas the
''straight lines'' ${\cal L}\left( Q_{0},{\bf P}_{0}{\bf P}_{1}\right) $,
drawn through the point $Q_{0}$ ($Q_{0}\neq P_{0}$) parallel to the vector $%
{\bf P}_{0}{\bf P}_{1}$, is not one-dimensional, in general. Note that the
Riemannian geometry denied a possibility of constructing the ''straight
line'' ${\cal L}\left( Q_{0},{\bf P}_{0}{\bf P}_{1}\right) $, referring to
lack of absolute parallelism. Lack of one-dimensionality for ${\cal L}\left(
Q_{0},{\bf P}_{0}{\bf P}_{1}\right) $ can be seen only in the $\sigma $%
-Riemannian geometry, which is defined as a consistent T-geometry, whose
world function is restricted by the relation (\ref{c1}). In the present
paper we consider and compare definitions of parallelism in Riemannian
geometry and in the consistent T-geometry ($\sigma $-Riemannian geometry)
and discuss corollaries of the Riemannian geometry inconsistency.

\section{Definition of parallelism}

Vector ${\bf P}_{0}{\bf P}_{1}\equiv \overrightarrow{P_{0}P_{1}}$ in
T-geometry is the ordered set of two points ${\bf P}_{0}{\bf P}_{1}=\left\{
P_{0},P_{1}\right\} $, \ \ $P_{0},P_{1}\in \Omega $. (The points $P_{0}$, $%
P_{1}$ may be similar). The scalar product $\left( {\bf P}_{0}{\bf P}_{1}.%
{\bf Q}_{0}{\bf Q}_{1}\right) $ of two vectors ${\bf P}_{0}{\bf P}_{1}$ and $%
{\bf Q}_{0}{\bf Q}_{1}$ is defined by the relation 
\begin{equation}
\left( {\bf P}_{0}{\bf P}_{1}.{\bf Q}_{0}{\bf Q}_{1}\right) =\sigma \left(
P_{0},Q_{1}\right) +\sigma \left( P_{1},Q_{0}\right) -\sigma \left(
P_{0},Q_{0}\right) -\sigma \left( P_{1},Q_{1}\right) ,  \label{a2}
\end{equation}
for all $P_{0},P_{1},Q_{0},Q_{1}\in \Omega $. As it follows from (\ref{a1}),
(\ref{a2}), in the symmetric T-geometry 
\begin{equation}
\left( {\bf P}_{0}{\bf P}_{1}.{\bf Q}_{0}{\bf Q}_{1}\right) =\left( {\bf Q}%
_{0}{\bf Q}_{1}.{\bf P}_{0}{\bf P}_{1}\right) ,\qquad \forall
P_{0},P_{1},Q_{0},Q_{1}\in \Omega  \label{a2.1}
\end{equation}
Further we shall consider only symmetric T-geometry and shall not stipulate
this. (Asymmetric T-geometry is considered in \cite{R002}).

When the world function $\sigma $ is such one \cite{R02} that the $\sigma $%
-space $V=\left\{ \sigma ,\Omega \right\} $ is the $n$-dimensional proper
Euclidean space $E_{n}$ the scalar product (\ref{a2}) turns to the scalar
product of two vectors in $E_{n}$. Besides, it follows from (\ref{a1}), (\ref
{a2}) that in any T-geometry 
\begin{equation}
\left( {\bf P}_{0}{\bf P}_{1}.{\bf Q}_{0}{\bf Q}_{1}\right) =-\left( {\bf P}%
_{1}{\bf P}_{0}.{\bf Q}_{0}{\bf Q}_{1}\right) ,\qquad \forall
P_{0},P_{1},Q_{0},Q_{1}\in \Omega  \label{a3}
\end{equation}
\begin{equation}
\left( {\bf P}_{0}{\bf P}_{1}.{\bf Q}_{0}{\bf Q}_{1}\right) +\left( {\bf P}%
_{1}{\bf P}_{2}.{\bf Q}_{0}{\bf Q}_{1}\right) =\left( {\bf P}_{0}{\bf P}_{2}.%
{\bf Q}_{0}{\bf Q}_{1}\right) ,  \label{a5.0}
\end{equation}
for all $P_{0},P_{1},P_{2},Q_{0},Q_{1}\in \Omega $. Two vectors ${\bf P}_{0}%
{\bf P}_{1}$ and ${\bf Q}_{0}{\bf Q}_{1}$ are parallel $\left( {\bf P}_{0}%
{\bf P}_{1}\uparrow \uparrow {\bf Q}_{0}{\bf Q}_{1}\right) $, if 
\begin{eqnarray}
\left( {\bf P}_{0}{\bf P}_{1}\uparrow \uparrow {\bf Q}_{0}{\bf Q}_{1}\right)
&:&\;\;\left( {\bf P}_{0}{\bf P}_{1}.{\bf Q}_{0}{\bf Q}_{1}\right) =\left| 
{\bf P}_{0}{\bf P}_{1}\right| \cdot \left| {\bf Q}_{0}{\bf Q}_{1}\right| ,
\label{a6} \\
\left| {\bf P}_{0}{\bf P}_{1}\right| &\equiv &\sqrt{\left( {\bf P}_{0}{\bf P}%
_{1}.{\bf P}_{0}{\bf P}_{1}\right) },\qquad \left| {\bf Q}_{0}{\bf Q}%
_{1}\right| \equiv \sqrt{\left( {\bf Q}_{0}{\bf Q}_{1}.{\bf Q}_{0}{\bf Q}%
_{1}\right) }  \nonumber
\end{eqnarray}
Definition of parallelism (\ref{a6}) does not contain a reference to
coordinate system, to a path of parallel transport, or to other means of
description. The relation (\ref{a6}) determines parallelism of two remote
vectors, using only world function $\sigma $. Parallelism of two vectors is
absolute in the sense that any two vectors ${\bf P}_{0}{\bf P}_{1}$ and $%
{\bf Q}_{0}{\bf Q}_{1}{\bf \ }$ are either parallel or not.

Vector ${\bf u}$ in $n$-dimensional Riemannian geometry is defined as a set
of $n$ quantities ${\bf u}=\left\{ u_{i}\right\} $,\ \ $i=1,2,...n,$ given
at some coordinate system $K_{n}$ with coordinates $x=\left\{ x^{i}\right\}
,\;\;i=1,2,...n$. At the coordinate transformation $K_{n}\rightarrow \tilde{K%
}_{n}$%
\begin{equation}
x^{i}\rightarrow \tilde{x}^{i}=\tilde{x}^{i}\left( x\right) ,\qquad
i=1,2,...n  \label{b1}
\end{equation}
covariant components $u_{i}$ of the vector ${\bf u}$ transforms as follows 
\begin{equation}
u_{i}\rightarrow \tilde{u}_{i}=\frac{\partial x^{k}}{\partial \tilde{x}^{i}}%
u_{k},\qquad i=1,2,...n  \label{b2}
\end{equation}
Summation from $1$ to $n$ is made over repeating indices.

Let $x$ be coordinates of the point $P$, and $x^{\prime }$ be coordinates of
the point $P^{\prime }$. Then the vector ${\bf PP}^{\prime }$ at the point $%
P $ is introduced by the relation 
\begin{equation}
{\bf PP}^{\prime }=\left\{ -\sigma _{i}\left( x,x^{\prime }\right) \right\}
,\qquad i=1,2,...n  \label{b3}
\end{equation}
\begin{equation}
\sigma _{i}\equiv \partial _{i}\sigma \left( x,x^{\prime }\right) \equiv 
\frac{\partial \sigma \left( x,x^{\prime }\right) }{\partial x^{i}},\qquad
i=1,2,...n  \label{b4}
\end{equation}
where the world function $\sigma $ is defined by the relation (\ref{c1}).
Here $\sigma \left( x,x^{\prime }\right) =\sigma \left( P,P^{\prime }\right) 
$ is the world function between the points $P$ and $P^{\prime }$.

In the $n$-dimensional proper Euclidean space $E_{n}$ and rectilinear
coordinate system $K_{n}$ the world function has the form 
\begin{equation}
\sigma \left( x,x^{\prime }\right) =\frac{1}{2}g_{\left( {\rm E}\right)
ik}\left( x^{i}-x^{\prime i}\right) \left( x^{k}-x^{\prime k}\right) ,\qquad
g_{\left( {\rm E}\right) ik}=\text{const}  \label{b7}
\end{equation}
and according to (\ref{b3}) the vector ${\bf PP}^{\prime }$ has covariant
coordinates $\left\{ g_{\left( {\rm E}\right) ik}\left( x^{k}-x^{\prime
k}\right) \right\} $, $\;i=1,2,...n$. Scalar product of two vectors ${\bf PP}%
^{\prime }$ and ${\bf PP}^{\prime \prime }$, having common origin at the
point $P$ has the form 
\begin{equation}
\left( {\bf PP}^{\prime }.{\bf PP}^{\prime \prime }\right) _{{\rm R}_{{\rm n}%
}}=g^{ik}\left( x\right) \sigma _{i}\left( x,x^{\prime }\right) \sigma
_{k}\left( x,x^{\prime \prime }\right)  \label{b9}
\end{equation}
where index ''R$_{{\rm n}}$'' means that the scalar product is defined in
the Riemannian space $R_{n}$ according to conventional rules of Riemannian
geometry.

According to (\ref{b9}) and in virtue of properties (\ref{c2}) of the world
function of the Riemannian space we obtain 
\begin{equation}
\left| {\bf PP}^{\prime }\right| ^{2}\equiv \left( {\bf PP}^{\prime }.{\bf PP%
}^{\prime }\right) =2\sigma \left( P,P^{\prime }\right)  \label{b6a}
\end{equation}
The definition (\ref{b9}) coincide with the general definition (\ref{a2}) in
the following cases: (1) if the Riemannian space $R_{n}$ coincide with the
Euclidean space $E_{n}$, (2) if vectors ${\bf PP}^{\prime }$ and ${\bf PP}%
^{\prime \prime }$ are infinitesimally small, (3) if $\sigma _{i}\left(
x,x^{\prime }\right) =a\sigma _{i}\left( x,x^{\prime \prime }\right) $, \ $%
i=1,2,...n,$ $a=$const (as it follows from (\ref{c2}), (\ref{b6a})). In
other cases the scalar products (\ref{b9}) and (\ref{a2}) do not coincide,
in general. Besides, the scalar product (\ref{b9}) is defined only for
vectors having a common origin. In the case of vectors ${\bf PP}^{\prime }$
and ${\bf QQ}^{\prime }$ with different origins the scalar product $\left( 
{\bf PP}^{\prime }.{\bf QQ}^{\prime }\right) $ must be defined in addition.
But this scalar product is not defined in Riemannian geometry, because to
define $\left( {\bf PP}^{\prime }.{\bf QQ}^{\prime }\right) $ for $Q\neq P$,
the vector ${\bf PP}^{\prime }$ must be transported at the point $Q$ in
parallel, and thereafter the definition (\ref{b9}) should be used. Result of
parallel transport depends on the path of transport, and the scalar product $%
\left( {\bf PP}^{\prime }.{\bf QQ}^{\prime }\right) $ for $Q\neq P$ cannot
be defined uniquely. If one uses definition (\ref{a2}) and relation (\ref{c1}%
) for determination of $\left( {\bf PP}^{\prime }.{\bf QQ}^{\prime }\right) $
for ${\bf Q}\neq {\bf P}$ the result is unique, but definition of
parallelism on the base of this scalar product leads to a set of many
vectors ${\bf QQ}^{\prime }$, which are parallel to ${\bf PP}^{\prime }$,
whereas the conventional conception of Riemannian geometry demands that such
a vector would be only one. In other words, the Riemannian geometry becomes
to be inconsistent at this point.

The definition (\ref{a2}) does not contain any reference to the means of
description, whereas the definition (\ref{b9}) does. The definition (\ref{b9}%
) is invariant with respect to coordinate transformation (\ref{b1}), but it
refers to the dimension $n$ of the space $R_{n}$ and existence of $n$%
-dimensional manifold. It means that the definition (\ref{a2}) is more
general and perfect, because it does not use special properties of the
Euclidean world function $\sigma _{{\rm E}}$.

These special properties of $n$-dimensional proper Euclidean space are
determined as follows \cite{R02}.

I: 
\begin{equation}
\exists {\cal P}^{n}=\left\{ P_{0},P_{1},...P_{n}\right\} ,\qquad
F_{n}\left( {\cal P}^{n}\right) \neq 0,\qquad F_{k}\left( {\Omega }%
^{k+1}\right) =0,\qquad k>n  \label{b10}
\end{equation}
where 
\begin{equation}
F_{n}\left( {\cal P}^{n}\right) =\det \left| \left| \left( {\bf P}_{0}{\bf P}%
_{i}.{\bf P}_{0}{\bf P}_{k}\right) \right| \right| =\det \left| \left|
g_{ik}\left( {\cal P}^{n}\right) \right| \right| \neq 0,\qquad i,k=1,2,...n
\label{b11}
\end{equation}
Vectors ${\bf P}_{0}{\bf P}_{i}$, $\;i=1,2,...n$ are basic vectors of the
rectilinear coordinate system $K_{n}$ with the origin at the point $P_{0}$,
and metric tensors $g_{ik}\left( {\cal P}^{n}\right) $, $g^{ik}\left( {\cal P%
}^{n}\right) $, \ $i,k=1,2,...n$ in $K_{n}$ are defined by relations 
\begin{equation}
g^{ik}\left( {\cal P}^{n}\right) g_{lk}\left( {\cal P}^{n}\right) =\delta
_{l}^{i},\qquad g_{il}\left( {\cal P}^{n}\right) =\left( {\bf P}_{0}{\bf P}%
_{i}.{\bf P}_{0}{\bf P}_{l}\right) ,\qquad i,l=1,2,...n  \label{a15b}
\end{equation}

II: 
\begin{equation}
\sigma _{{\rm E}}\left( P,Q\right) =\frac{1}{2}g^{ik}\left( {\cal P}%
^{n}\right) \left( x_{i}\left( P\right) -x_{i}\left( Q\right) \right) \left(
x_{k}\left( P\right) -x_{k}\left( Q\right) \right) ,\qquad \forall P,Q\in 
{\Bbb R}^{n}  \label{a15a}
\end{equation}
where coordinates $x_{i}\left( P\right) $ of the point $P$ are defined by
the relation 
\begin{equation}
x_{i}\left( P\right) =\left( {\bf P}_{0}{\bf P}_{i}.{\bf P}_{0}{\bf P}%
\right) ,\qquad i=1,2,...n  \label{b12}
\end{equation}

III: The metric tensor matrix $g_{lk}\left( {\cal P}^{n}\right) $ has only
positive eigenvalues 
\begin{equation}
g_{k}>0,\qquad k=1,2,...,n  \label{a15c}
\end{equation}

IV. Continuity condition: the system of equations 
\begin{equation}
\left( {\bf P}_{0}{\bf P}_{i}.{\bf P}_{0}{\bf P}\right) =y_{i}\in {\Bbb R}%
,\qquad i=1,2,...n  \label{b.12}
\end{equation}
considered to be equations for determination of the point $P$ as a function
of coordinates $y=\left\{ y_{i}\right\} $,\ \ $i=1,2,...n$ has always one
and only one solution.

Conditions I -- IV are necessary and sufficient conditions of that the $%
\sigma $-space $V=\left\{ \sigma ,\Omega \right\} $ is the $n$-dimensional
proper Euclidean space \cite{R02}. These special properties of $E_{n}$ are
different for different dimension $n$, and contain a reference to $n$.

Let us use in Riemannian geometry two different definitions of parallelism,
based on application of relations (\ref{a6}), (\ref{a2}) and (\ref{a6}), (%
\ref{b9}) respectively. Although definitions of (\ref{a2}) and (\ref{b9})
for the scalar product are different, they give the same result for
parallelism of to vectors having a common origin.

The relations (\ref{a6}), (\ref{b9}) define, parallelism only for two
vectors, having a common origin. To define parallelism of two remote vectors 
${\bf u}\left( x\right) $ and ${\bf u}\left( x^{\prime }\right) $ in
Riemannian geometry, one defines parallelism of two infinitesimally close
vectors ${\bf u}\left( x\right) $ and ${\bf u}\left( x+dx\right) $ by means
of the relation 
\begin{equation}
u_{i}\left( x+dx\right) =u_{i}\left( x\right) -\Gamma _{il}^{k}\left(
x\right) u_{k}\left( x\right) dx^{l},\qquad i=1,2,...n  \label{b41}
\end{equation}
\begin{equation}
\Gamma _{il}^{k}=\frac{1}{2}g^{kj}\left( g_{ij,l}+g_{lj,i}-g_{il,j}\right)
,\qquad g_{ij,l}\equiv \frac{\partial g_{ij}}{\partial x^{l}}  \label{b42}
\end{equation}
The vector ${\bf u}\left( x^{\prime }\right) $ at the point $x^{\prime }$
parallel to the vector ${\bf u}\left( x\right) $ is obtained by subsequent
application of the infinitesimally small transport (\ref{b41}) along some
path ${\cal L}$, connecting points $x$ and $x^{\prime }$. Note that the
vectors ${\bf u}\left( x\right) $ and ${\bf u}\left( x+dx\right) $ are
parallel, and besides they have the same length. In general, result of the
parallel transport along ${\cal L}$ depends on ${\cal L}$. Such a situation
is known as a lack of absolute parallelism. For flat Riemannian spaces there
is the absolute parallelism, but for the curved Riemannian spaces the
absolute parallelism is absent, in general.

Application of the parallelism definition, based on relations (\ref{a6}), (%
\ref{a2}), to vectors ${\bf PP}^{\prime }$ and ${\bf P}_{1}{\bf P}^{\prime
\prime }$ in Riemannian geometry with infinitesimally close points $P$ and $%
P_{1}$ gives a result coinciding with (\ref{b41}), only if the displacement
vector ${\bf PP}_{1}||{\bf PP}^{\prime }$ (and hence ${\bf PP}_{1}||{\bf PP}%
^{\prime \prime }$). This property provides one-dimensionality of geodesics,
obtained as a result of deformation of Euclidean straight lines. In other
cases, the results of two definitions of parallelism appear to be different,
in general, because the relation (\ref{b41}) gives only one vector ${\bf u}%
\left( x+dx\right) ,$ parallel to ${\bf u}\left( x\right) $, whereas
relations (\ref{a6}), (\ref{a2}) generate, in general, a set of many vectors 
${\bf P}_{1}{\bf P}^{\prime \prime }$, which are parallel to ${\bf PP}%
^{\prime }$, but which are not parallel, in general, between themselves \cite
{R02}. The difference is conditioned by the fact that the condition of
parallelism (\ref{a6}) contains only one relation, whereas the condition of
parallelism (\ref{b41}) contains $n$ relations.

To explain the reason of this difference, let us consider the case, when $%
\left| {\bf PP}^{\prime }\right| $ $\neq 0$ and $\left| {\bf PP}^{\prime
\prime }\right| \neq 0$. In this case one can itroduce unit vectors \quad $%
\sigma _{i}\left( x,x^{\prime }\right) \left( 2\sigma \left( x,x^{\prime
}\right) \right) ^{-1/2}$, \ $\sigma _{i}\left( x,x^{\prime \prime }\right)
\left( 2\sigma \left( x,x^{\prime \prime }\right) \right) ^{-1/2}$ and
rewrite relations (\ref{a6}), (\ref{b9}) in the form of scalar product of
the two unit vectors 
\begin{equation}
g^{ik}\left( x\right) \frac{\sigma _{i}\left( x,x^{\prime }\right) }{\sqrt{%
2\sigma \left( x,x^{\prime }\right) }}\frac{\sigma _{k}\left( x,x^{\prime
\prime }\right) }{\sqrt{2\sigma \left( x,x^{\prime \prime }\right) }}=1,
\label{b43}
\end{equation}
Let the matrix of metric tensor $g^{ik}\left( x\right) $ has eigenvalues of
the same sign. Then both vectors $\sigma _{i}\left( x,x^{\prime }\right)
\left( 2\sigma \left( x,x^{\prime }\right) \right) ^{-1/2}$ and $\sigma
_{i}\left( x,x^{\prime \prime }\right) \left( 2\sigma \left( x,x^{\prime
\prime }\right) \right) ^{-1/2}$ are equal, and one relation (\ref{b43}) is
equivalent to $n$ relations 
\begin{equation}
\sigma _{i}\left( x,x^{\prime }\right) =a\sigma _{i}\left( x,x^{\prime
\prime }\right) ,\qquad i=1,2,...n,\qquad a>0  \label{a9}
\end{equation}
where $a$ is some constant. Conditions (\ref{a9}) with arbitrary $a\neq 0$
mean that vectors ${\bf PP}^{\prime }$ and ${\bf PP}^{\prime \prime }$,
having a common origin, are collinear (parallel or antiparallel), provided
their components are proportional.

In the $n$-dimensional proper Euclidean space $E_{n}$ this condition can be
written $\sigma $-immanently. Let vector ${\bf P}_{0}{\bf R}$ be collinear
to the vector ${\bf P}_{0}{\bf P}_{1}$. Let us choose $n-1$ points $\left\{
P_{2},P_{3},...P_{n}\right\} $ in such a way, that $n$ vectors ${\bf P}_{0}%
{\bf P}_{i}$,$\;\;i=1,2,...n$ form a basis. Then the collinearity condition (%
\ref{a9}) of vectors ${\bf P}_{0}{\bf R}$ and ${\bf P}_{0}{\bf P}_{1}$ takes
the form of $n$ relations 
\begin{equation}
\left( {\bf P}_{0}{\bf P}_{i}.{\bf P}_{0}{\bf R}\right) =a\left( {\bf P}_{0}%
{\bf P}_{i}.{\bf P}_{0}{\bf P}_{1}\right) ,\qquad i=1,2,...n  \label{a10}
\end{equation}
Eliminating $a$ from $n$ relations (\ref{a10}) we obtain $n-1$ relations,
which are written in the form 
\begin{equation}
{\bf P}_{0}{\bf P}_{1}||{\bf P}_{0}{\bf R}:\qquad \left| 
\begin{array}{cc}
\left( {\bf P}_{0}{\bf P}_{1}.{\bf P}_{0}{\bf R}\right) & \left( {\bf P}_{0}%
{\bf P}_{i}.{\bf P}_{0}{\bf R}\right) \\ 
\left( {\bf P}_{0}{\bf P}_{1}.{\bf P}_{0}{\bf P}_{1}\right) & \left( {\bf P}%
_{0}{\bf P}_{i}.{\bf P}_{0}{\bf P}_{1}\right)
\end{array}
\right| =0,\qquad i=2,3,...n  \label{a12}
\end{equation}
Thus, we have two different formulation of the collinearity conditions of
vectors ${\bf P}_{0}{\bf R}$ and ${\bf P}_{0}{\bf P}_{1}$: (\ref{a12}) and
the relation 
\begin{equation}
{\bf P}_{0}{\bf P}_{1}||{\bf P}_{0}{\bf R}:\qquad \left| 
\begin{array}{cc}
\left( {\bf P}_{0}{\bf P}_{1}.{\bf P}_{0}{\bf P}_{1}\right) & \left( {\bf P}%
_{0}{\bf P}_{1}.{\bf P}_{0}{\bf R}\right) \\ 
\left( {\bf P}_{0}{\bf R}.{\bf P}_{0}{\bf P}_{1}\right) & \left( {\bf P}_{0}%
{\bf R}.{\bf P}_{0}{\bf R}\right)
\end{array}
\right| =0  \label{b44}
\end{equation}
which follows from (\ref{a6}). In $E_{n}$ conditions (\ref{a12}), and (\ref
{b44}) are equivalent, because the choice of $n-1$ points $\left\{
P_{2},P_{3},...P_{n}\right\} $ is arbitrary, and they are fictitious in (\ref
{a12}). The collinearity conditions (\ref{a12}) and (\ref{b44}) are
equivalent due to special properties (\ref{a15a}) of $E_{n}$. In the $n$%
-dimensional proper Riemannian geometry the conditions (\ref{a12}), and (\ref
{b44}) are also equivalent, and points $\left\{ P_{2},P_{3},...P_{n}\right\} 
$ are also fictitious in (\ref{a12}). This is connected with the special
choice of the world function (\ref{c1}) of $n$-dimensional Riemannian space.
At another choice of the world function the points $\left\{
P_{2},P_{3},...P_{n}\right\} $ stop to be fictitious.

To manifest difference between the conditions (\ref{a12}) and (\ref{b44}),
let us construct the ''straight line'' ${\cal T}_{P_{0}P_{1}}$, passing
through points $P_{0},P_{1}$, defining it as set of such points $R$, that $%
{\bf P}_{0}{\bf R}||{\bf P}_{0}{\bf P}_{1}$. Using two variants of the
collinearity conditions (\ref{a12}), and (\ref{b44}) we obtain two different
geometric objects 
\begin{equation}
{\cal T}_{P_{0}P_{1}}=\left\{ R|\;{\bf P}_{0}{\bf P}_{1}||{\bf P}_{0}{\bf R}%
\right\} =\left\{ R|\left( {\bf P}_{0}{\bf P}_{1}.{\bf P}_{0}{\bf R}\right)
^{2}=\left| {\bf P}_{0}{\bf P}_{1}\right| ^{2}\left| {\bf P}_{0}{\bf R}%
\right| ^{2}\right\}  \label{a8}
\end{equation}
and 
\begin{equation}
{\cal L}=\left\{ R\left| \bigwedge\limits_{k=2}^{k=n}f\left(
P_{0},P_{1},P_{k},R\right) =0\right. \right\}
=\bigcap\limits_{k=2}^{k=n}\left\{ R\left| f\left(
P_{0},P_{1},P_{k},R\right) =0\right. \right\}  \label{a14}
\end{equation}
where 
\begin{equation}
f\left( P_{0},P_{1},P_{i},R\right) =\left| 
\begin{array}{cc}
\left( {\bf P}_{0}{\bf P}_{1}.{\bf P}_{0}{\bf R}\right) & \left( {\bf P}_{0}%
{\bf P}_{i}.{\bf P}_{0}{\bf R}\right) \\ 
\left( {\bf P}_{0}{\bf P}_{1}.{\bf P}_{0}{\bf P}_{1}\right) & \left( {\bf P}%
_{0}{\bf P}_{i}.{\bf P}_{0}{\bf P}_{1}\right)
\end{array}
\right| =0,\qquad i=2,3,...n  \label{b45}
\end{equation}
In the $n$-dimensional proper Euclidean space and in the $n$-dimensional
proper Riemannian space the geometric objects ${\cal L}$ and ${\cal T}%
_{P_{0}P_{1}}$ coincide, but at a more general form of the world function
the geometric objects ${\cal L}$ and ${\cal T}_{P_{0}P_{1}}$ are different,
in general.

The relation (\ref{a14}) determines the straight line ${\cal L}$ in the $n$%
-dimensional proper Euclidean space as an intersection of $n-1$ $\left(
n-1\right) $-dimensional surfaces 
\begin{equation}
{\cal S}\left( P_{0},P_{1},P_{k}\right) =\left\{ R\left| f\left(
P_{0},P_{1},P_{k},R\right) =0\right. \right\} ,\qquad k=2,3,...n  \label{a15}
\end{equation}
In general, such an intersection is a one-dimensional line, but this line is
determined by $n+1$ points ${\cal P}^{n}\equiv \left\{
P_{0},P_{1},...,P_{n}\right\} $, whereas the ''straight line'' ${\cal T}%
_{P_{0}P_{1}}$, defined by the relation (\ref{a8}), depends only on two
points $P_{0},P_{1}$.

In general case, when the special properties of the Euclidean space
disappear, the relation (\ref{a14}) describes one-dimensional object
depending on more than two points. Thus, one can eliminate dependence of the
collinearity definition (\ref{a9}) on the coordinate system, but instead of
this dependence a dependence on additional points appears. These additional
points $P_{2},P_{3},...$ represent the coordinate system in the $\sigma $%
-immanent form. The number of additional points which are necessary for
determination of the ''straight line'' (\ref{a9}) as a one-dimensional line
depends on the dimension of the Euclidean space. From formal viewpoint the
geometric object ${\cal L}$, determined $\sigma $-immanently by (\ref{a14}),
is not a straight line, but some other geometric object, coinciding with the
straight line in the $n$-dimensional proper Euclidean space.

The straight line in the $n$-dimensional proper Euclidean space has two
properties: (1) the straight line is determined by two points $P_{0}$, $%
P_{1} $ independently of the dimension of the Euclidean space, (2) the
straight line is a one-dimensional line. In general, both properties are not
retained at deformation of the Euclidean space. If we use the definition (%
\ref{a8}), we retain the first property, but violate, in general, the second
one. If we use the definition (\ref{a14}), depending on the Euclidean space
dimension and on the way of description (in the form of coordinate system,
or in the form of additional arbitrary points), we retain the second
property and violate, in general, the first one. Which of the two
definitions of the ''straight line'' should be used?

The answer is evident. Firstly, the definition (\ref{a8}) does not refer to
any means of description, whereas the definition (\ref{a14}) does. Secondly,
the property of the ''straight line'' of being determined by two points is
the more natural property of geometry, than the property of being a
one-dimensional line. Use of the definition (\ref{a8}) is a logical
necessity, but not a hypothesis, which can be confirmed or rejected in
experiment. Consideration of the ''straight line'' as a one-dimensional
geometric object {\it in any geometry} is simply a preconception, based on
the fact, that in the proper Euclidean geometry the straight line is a
one-dimensional geometric object. The statement that there is only one
vector ${\bf Q}_{0}{\bf Q}_{1}$ of fixed length which is parallel to the
vector ${\bf P}_{0}{\bf P}_{1}$ is another formulation of the preconception
mentioned above.

\section{Consequence of inconsistent definition \newline
of parallelism}

Abstracting from the history of the Riemannian geometry creation and motives
of its creation, let us evaluate what is the Riemannian geometry as a kind
of physical geometry. The conventional Riemannian geometry is to be a
special case of a physical geometry, constructed on the basis of the
principle of geometry deformation. The Riemannian geometry uses definition
of the scalar product (\ref{b9}), which is completely compatible with the
principle of geometry deformation only for several geometries. To compensate
inconsistencies, generated by incorrectness of definition (\ref{b9}), the
Riemannian geometry uses the constraint (\ref{c1}), tending to eliminate
geometries, for which the definition (\ref{b9}) is inconsistent. The
constraint (\ref{c1}) removes most of possible inconsistencies, but not all,
and the Riemannian geometry appears to be inconsistent geometry.

In the contemporary geometry and physics the definition (\ref{a9}) or (\ref
{a14}) is used, and this circumstance is a reason for many problems, because
this definition lies in the foundation of the geometry, and the geometry in
turn lies in the foundation of physics.

Let us list some consequences of the statement that the straight line is a
one-dimensional geometric object in any space-time geometry.

\begin{enumerate}
\item  Lack of absolute parallelism in the space-time geometry (i.e. in
Riemannian geometry used for description of the space-time).

\item  Discrimination of any space-time geometry, where the timelike
straight is not a one-dimensional object, and (as a corollary)
discrimination of stochastic motion of microparticles.

\item  Consideration of spacelike straights, describing superlight particles
(tachyons), in the Minkowski space-time geometry, as one-dimensional
geometric objects.
\end{enumerate}

Let us discuss the first point. The world function of the Riemannian
geometry is chosen in such a way that the tube ${\cal T}_{P_{0}P_{1}}$ (we
use this term instead of the term ''straight line''), passing through the
points $P_{0}$, $P_{1}$ and defined by the relation (\ref{a8}), is a
one-dimensional geometric object in the Riemannian space-time geometry,
provided interval between the points $P_{0}$, $P_{1}$ is timelike $\left(
\sigma \left( P_{0},P_{1}\right) >0\right) $. But the timelike tube 
\begin{equation}
{\cal T}\left( P_{0},P_{1};Q_{0}\right) =\left\{ R|\;{\bf P}_{0}{\bf P}_{1}||%
{\bf Q}_{0}{\bf R}\right\} =\left\{ R|\left( {\bf P}_{0}{\bf P}_{1}.{\bf Q}%
_{0}{\bf R}\right) ^{2}=\left| {\bf P}_{0}{\bf P}_{1}\right| ^{2}\left| {\bf %
Q}_{0}{\bf R}\right| ^{2}\right\}  \label{a16}
\end{equation}
passing through the point $Q_{0}$ parallel to the remote timelike vector $%
{\bf P}_{0}{\bf P}_{1}$, is not a one-dimensional object, in general, in the 
$\sigma $-Riemannian geometry (in Riemannian geometry ${\cal T}\left(
P_{0},P_{1};Q_{0}\right) $ is not defined). One cannot achieve that any
timelike tube (\ref{a8}) to be a one-dimensional geometric object. In other
words, one cannot suppress globally nondegeneracy of all collinearity cones
of timelike vectors ${\bf Q}_{0}{\bf R}$, parallel to the timelike vector $%
{\bf P}_{0}{\bf P}_{1}$, although locally the collinearity cone
nondegeneracy of timelike vectors ${\bf P}_{0}{\bf R}$, parallel to the
timelike vector ${\bf P}_{0}{\bf P}_{1}$, can be suppressed, if the world
function is restricted by the constraint (\ref{c1}). In fact, according to
the correct definition (\ref{b44}) in the $\sigma $-Riemannian geometry
there are many timelike vectors ${\bf Q}_{0}{\bf R}$ of fixed length, which
are parallel to the remote timelike vector ${\bf P}_{0}{\bf P}_{1}$. As far
as according to the Riemannian conception of geometry there is to be only
one timelike vector ${\bf Q}_{0}{\bf R}$ of fixed length, which is parallel
to the remote timelike vector ${\bf P}_{0}{\bf P}_{1}$, one cannot choose
one vector among the set of equivalent vectors ${\bf Q}_{0}{\bf R}$, and one
is forced to deny the absolute parallelism.

The point two. The Minkowski space-time geometry $T_{{\rm M}}$ with the $%
\sigma $-space $\left\{ \sigma _{{\rm M}},{\Bbb R}^{4}\right\} $ is the
unique uniform isotropic flat geometry in the class of Riemannian
geometries. The class of uniform isotropic T-geometries on the set ${\Bbb R}%
^{4}$ of points is described by the world function $\sigma =\sigma _{{\rm M}%
}+D\left( \sigma _{{\rm M}}\right) $, where the arbitrary distortion
function $D$ describes character of nondegeneracy of timelike tubes ${\cal T}%
_{P_{0}P_{1}}$. In the Minkowski space-time geometry a motion of free
particles is deterministic. If $D>0$ the world line of a free particle
appears to be stochastic, because the running point moves along the world
line in the direction of vector tangent to the world line. There are many
vectors tangent to the world line. The particle can move along any of them,
and its motion becomes stochastic, (see details in \cite{R91}). In fact,
motion of microparticles (electrons, protons, etc.) is stochastic. It means
that the Minkowski geometry is not a true space-time geometry. One should
choose such a space-time geometry, which could explain stochastic motion of
microparticles. Such a space-time geometry is possible. In this space-time
geometry the distortion function $D\left( \sigma _{{\rm M}}\right) =\hbar
/\left( 2bc\right) $ for $\sigma _{{\rm M}}>\sigma _{0}\approx \hbar /\left(
2bc\right) $, where $\hbar $ is the quantum constant, $c$ is the speed of
the light, and $b$ is a new universal constant. In such a space-time
geometry the world function contains the quantum constant $\hbar $, and
nonrelativistic quantum effects are explained as geometric effects \cite{R91}%
. Insisting on the definition (\ref{a9}) of the parallelism, we discriminate
space-time geometries with $D\neq 0$. As a result we are forced to use
incorrect space-time geometry and to explain quantum effects by additional
hypotheses (quantum principles).

Let us consider ''straight lines'' in the Minkowski geometry. Let us define
the ''straight line'' by the relation (\ref{a8}). Let ${\bf e=P}_{0}{\bf P}%
_{1}$ and ${\bf x}={\bf P}_{0}{\bf R}$ be the running vector. Then the
relation determining the ''straight line'' ${\cal T}_{P_{0}P_{1}}$ has the
form 
\begin{equation}
{\cal T}_{P_{0}P_{1}}:\qquad \left| 
\begin{array}{cc}
\left( {\bf e.e}\right) & \left( {\bf e}.{\bf x}\right) \\ 
\left( {\bf x}.{\bf e}\right) & \left( {\bf x}.{\bf x}\right)
\end{array}
\right| =0  \label{a17}
\end{equation}
Looking for its solution in the form 
\begin{equation}
{\bf x}={\bf e}\tau +{\bf y}  \label{a17.a}
\end{equation}
and substituting this expression in (\ref{a17}), we obtain the equation of
the same form. 
\begin{equation}
\left| 
\begin{array}{cc}
\left( {\bf e.e}\right) & \left( {\bf e}.{\bf y}\right) \\ 
\left( {\bf y}.{\bf e}\right) & \left( {\bf y}.{\bf y}\right)
\end{array}
\right| =0  \label{a18}
\end{equation}
Evident solution ${\bf y=}\alpha {\bf e}$ is not interesting, because it has
been taken into account in (\ref{a17.a}). Imposing constraint $\left( {\bf e}%
.{\bf y}\right) =0$, one obtains from (\ref{a18}) 
\[
\left( {\bf e}.{\bf y}\right) =0,\qquad {\bf y}^{2}=0 
\]
If the vector ${\bf e}$ is timelike, for instance, ${\bf e=}\left\{
1,0,0,0\right\} $, then ${\bf y}=0$. If the vector ${\bf e}$ is spacelike,
for instance$,{\bf e=}\left\{ 0,1,0,0\right\} $, then the solution has the
form ${\bf y}$ =$\left\{ a,0,a\cos \psi ,a\sin \psi \right\} $, where $a$
and $\psi $ are arbitrary parameters. Thus, in the Minkowski space the
timelike ''straight line'' is a one-dimensional object, whereas the
spacelike ''straight line''\ is a three-dimensional surface, containing the
one-dimensional spacelike straight line ${\bf x}={\bf e}\tau $. In other
words, timelike directions are degenerate, and free particles, moving with
the speed $v<c$, are described by one-dimensional timelike ''straight
lines''. The spacelike directions are nondegenerate, and free particles,
moving with the speed $v>c$ (tachyons) are described by three-dimensional
surfaces. It is difficult to say, what does it mean practically. But, maybe,
tachyons were not discovered, because they were searched in the form of
one-dimensional spacelike lines.

\end{document}